\documentclass[12pt,thmsa]{article}
\usepackage{amsmath,amssymb,amsthm}

\newtheorem{theorem}{Theorem}
\newtheorem{proposition}{Proposition}
\newtheorem{example}[theorem]{Example}

\begin{document}

\title{Multivector Functionals\thanks{%
published: \emph{Advances in Applied Clifford Algebras} \textbf{11}(S3),
93-103 (2001).}}
\author{A. M. Moya$^{1}\thanks{%
e-mail: moya@ime.unicamp.br}$, V. V. Fern\'{a}ndez$^{1}\thanks{%
e-mail: vvf@ime.unicamp.br}$ and W. A. Rodrigues Jr.$^{1,2}\thanks{%
e-mail: walrod@ime.unicamp.br or walrod@mpc.com.br}$ \\
$\hspace{-0.5cm}^{1}$ Institute of Mathematics, Statistics and Scientific
Computation\\
IMECC-UNICAMP CP 6065\\
13083-970 Campinas-SP, Brazil\\
$^{2}$ Department of Mathematical Sciences, University of Liverpool\\
Liverpool, L69 3BX, UK}
\date{11/26/2001 }
\maketitle

\begin{abstract}
In this paper we introduce the concept of \emph{multivector functionals.} We
study some possible kinds of derivative operators that can act in
interesting ways on these objects such as, e.g., the $A$-directional
derivative and the generalized concepts of curl, divergence and gradient.
The derivation rules are rigorously proved. Since the subject of this paper
has not been developed in previous literature, we work out in details
several examples of derivation of multivector functionals.
\end{abstract}

\tableofcontents

\section{Introduction}

This is the last paper (VII) of series of papers dealing with the theory of
multivector and extensor functions and multivector functionals. It is
dedicated to the introduction of a key concept, that of \emph{multivector
functionals} and the study of their properties. Particularly important is
the concept of \emph{induced multivector functionals.} Several kinds of
derivatives of multivector functionals, such as $A$-directional derivative
and generalized concepts of curl, divergence and gradient are defined. Since
the subject of the present paper has not been explored in the literature%
\footnote{%
For the best of our knowldge the only place where the concept has been
rudimentary used was in \cite{1}. The concept has been used also in (\cite{2}%
,\cite{3}).}, we present in section 3 several examples worked in detail of
calculations of different types of derivatives for multivector functionals.
Multivector functionals are fundamental for the formulation of the
Lagrangian field theory of multivector and extensor fields on an arbitrary
manifold, a subject that will be studied in a new series of papers.

\section{Multivector Functionals}

Any mapping which sends general extensors over $V$ into multivectors over $V$
will be called a \emph{general multivector functional over }$V$\emph{.}

In particular, the general functionals with image-values belonging to $%
\bigwedge^{r}V$ are said to be $r$\emph{-vector functionals of general
extensor.} For the cases $r=0,$ $r=1,$ $r=2,\ldots$ and $r=n$ we speak about
\emph{scalar, vector, bivector,}$\ldots$\emph{\ }and\emph{\ pseudoscalar
functionals}, respectively.

For the applications we have in mind we shall need only some particular
cases of these general functionals for which we will give special names.

Any mapping $\mathcal{F}:ext_{p}^{q}(V)\rightarrow\bigwedge^{r}V$ will be
called a $r$\emph{-vector functional of a }$(p,q)$\emph{-extensor.} In
accordance to what was said above, the cases for which $\mathcal{F}[t]$
belongs to $\mathbb{R},$ $V,$ $\bigwedge^{2}V,\ldots$ and $\bigwedge^{n}V$
will be named respectively as \emph{scalar, vector, bivector,}$\ldots$ and
\emph{pseudoscalar functionals of a }$(p,q)$\emph{-extensor.}

\subsection{Induced Multivector Functionals}

Let $F:\underset{k\text{ factors}}{\underbrace{\bigwedge^{q}V\times
\cdots\times\bigwedge^{q}V}}\rightarrow\bigwedge^{r}V$ be any $r$-vector
function of $k$ $q$-vector variables. Take some $k$-uple of $p$-vectors $%
(A^{1},\ldots,A^{k}).$

Associated to $F$ and with respect to $(A^{1},\ldots,A^{k})$ it is possible
to construct a $r$-vector functional of a $(p,q)$-extensor, say $\mathcal{F}%
_{(A^{1},\ldots,A^{k})},$ given by
\begin{equation*}
ext_{p}^{q}(V)\ni t\mapsto\mathcal{F}_{(A^{1},\ldots,A^{k})}[t]\in
\bigwedge\nolimits^{r}V\text{ such that}
\end{equation*}

\begin{equation}
\mathcal{F}_{(A^{1},\ldots,A^{k})}[t]=F[t(A^{1}),\ldots,t(A^{k})].
\label{7.2a}
\end{equation}
It will be called the $r$\emph{-vector functional of a }$(p,q)$\emph{%
-extensor induced by }$F,$\emph{\ relative to }$(A^{1},\ldots,A^{k})$\emph{.}

If $F$ is differentiable on $\underset{k\text{ factors}}{\underbrace {%
\bigwedge\nolimits^{q}V\times\cdots\times\bigwedge^{q}V}},$ then $\mathcal{F}%
_{(A^{1},\ldots,A^{k})}$ is said to be \emph{differentially-induced by }$F$%
\emph{\ with respect to }$(A^{1},\ldots,A^{k}).$

In this way if $\mathcal{F}_{(A^{1},\ldots,A^{k})}$ is
differentially-induced, then there must exist the standard derivatives of $F$
with respect to each $p$-vector variable $X^{1},\ldots$ and $X^{k}$ (the
so-called \emph{partial derivatives of }$F$), i.e., $\partial_{X^{1}}F,\ldots
$ and $\partial_{X^{k}}F,$ (see \cite{6}).

Associated to $\partial_{X^{1}}F,\ldots$ and $\partial_{X^{k}}F$ with
respect to $(A^{1},\ldots,A^{k})$ we can define the following multivector
functionals of a $(p,q)$-extensor:
\begin{align*}
ext_{p}^{q}(V)\ni t &
\mapsto\partial_{X^{1}}F[t(A^{1}),\ldots,t(A^{k})]\in\bigwedge V, \\
& \ldots
\end{align*}
and
\begin{equation}
ext_{p}^{q}(V)\ni
t\mapsto\partial_{X^{k}}F[t(A^{1}),\ldots,t(A^{k})]\in\bigwedge V.
\label{7.2b}
\end{equation}
We see that they are induced by the partial derivatives of $F$ with respect
to $(A^{1},\ldots,A^{k}).$

\subsubsection{Directional Derivative}

Take an arbitrary $p$-vector $A.$ We introduce the\emph{\ }$A$\emph{%
-directional derivative} of the differentially-induced $r$-vector functional
$\mathcal{F}_{(A^{1},\ldots,A^{k})}$ as being the multivector functional $%
\mathcal{F}_{(A^{1},\ldots,A^{k})A}^{\prime}$ given by
\begin{equation*}
ext_{p}^{q}(V)\ni t\mapsto\mathcal{F}_{(A^{1},\ldots,A^{k})A}^{\prime}[t]\in%
\Lambda V\text{ such that}
\end{equation*}

\begin{equation}
\mathcal{F}_{(A^{1},\ldots,A^{k})A}^{\prime}[t]=\overset{k}{\underset{i=1}{%
\sum}}A\cdot A^{i}\partial_{X^{i}}F[t(A^{1}),\ldots,t(A^{k})].   \label{7.3a}
\end{equation}

Note that the algebraic object just defined associated to $\mathcal{F}%
_{(A^{1},\ldots,A^{k})}$ has the property of linearity with respect to the
direction, i.e., for any $\alpha,\beta\in\mathbb{R}$ and $A,B\in\bigwedge
^{p}V$
\begin{equation}
\mathcal{F}_{(A^{1},\ldots,A^{k})\alpha A+\beta B}^{\prime}[t]=\alpha
\mathcal{F}_{(A^{1},\ldots,A^{k})A}^{\prime}[t]+\beta\mathcal{F}%
_{(A^{1},\ldots,A^{k})B}^{\prime}[t],   \label{7.3c}
\end{equation}
as expected to hold for a well-defined $A$-directional derivative of $%
\mathcal{F}_{(A^{1},\ldots,A^{k})}.$

\subsubsection{Derivatives}

Let $(\{e_{k}\},\{e^{k}\})$ be a pair of arbitrary reciprocal bases of $V.$
It is also possible to introduce \emph{four} derivatives-like operators for
the differentially-induced $r$-vector functional $\mathcal{F}_{(A^{1},\ldots
,A^{k})}$ as the following multivector functionals $*\mathcal{F}%
_{(A^{1},\ldots,A^{k})}^{\prime}$ defined by
\begin{align}
\ast\mathcal{F}_{(A^{1},\ldots,A^{k})}^{\prime}[t] & =\frac{1}{p!}%
(e^{j_{1}}\wedge\ldots e^{j_{p}})*\mathcal{F}_{(A^{1},\ldots,A^{k})e_{j_{1}}%
\wedge\ldots e_{j_{p}}}^{\prime}[t]  \label{7.4aa} \\
& =\frac{1}{p!}(e_{j_{1}}\wedge\ldots e_{j_{p}})*\mathcal{F}%
_{(A^{1},\ldots,A^{k})e^{j_{1}}\wedge\ldots e^{j_{p}}}^{\prime}[t],
\label{7.4a}
\end{align}
where $*$ means either $(\wedge),$ $(\cdot),$ $(\lrcorner)$ or $($\emph{%
Clifford product}$).$

It should be noted that $*\mathcal{F}_{(A^{1},\ldots,A^{k})}^{\prime}$ are
well-defined multivector functionals of $(p,q)$-extensor \emph{only }%
associated with $\mathcal{F}_{(A^{1},\ldots,A^{k})}$ since, by taking into
account eq.(\ref{7.3c}), $*\mathcal{F}_{(A^{1},\ldots,A^{k})}^{\prime}[t]$
are multivectors which do not depend on the choice of $(\{e_{k}\},\{e^{k}\})$%
.

Recall also that a straightforward calculation gives with the use eq.(\ref%
{7.3a}) that
\begin{align}
\ast\mathcal{F}_{(A^{1},\ldots,A^{k})}^{\prime}[t] & =\frac{1}{p!}%
(e^{j_{1}}\wedge\ldots e^{j_{p}})*(\overset{k}{\underset{i=1}{\sum}}%
(e_{j_{1}}\wedge\ldots e_{j_{p}})\cdot A^{i}\partial_{X^{i}}F[\ldots])
\notag \\
& =(\overset{k}{\underset{i=1}{\sum}}\frac{1}{p!}(e_{j_{1}}\wedge\ldots
e_{j_{p}})\cdot A^{i}e^{j_{1}}\wedge\ldots
e^{j_{p}})*\partial_{X^{i}}F[\ldots]  \notag \\
\ast\mathcal{F}_{(A^{1},\ldots,A^{k})}^{\prime}[t] & =\overset{k}{\underset{%
i=1}{\sum}}A^{i}*\partial_{X^{i}}F[t(A^{1}),\ldots,t(A^{k})].   \label{7.4b}
\end{align}
Eq.(\ref{7.4b}) shows explicitly that $*\mathcal{F}_{(A^{1},\ldots,A^{k})}^{%
\prime}$ can be intrinsically defined without using any pair of reciprocal
bases of $V$.

The special cases: $\wedge\mathcal{F}_{(A^{1},\ldots,A^{k})}^{\prime},$ $%
\cdot\mathcal{F}_{(A^{1},\ldots,A^{k})}^{\prime},$ $\lrcorner\mathcal{F}%
_{(A^{1},\ldots,A^{k})}^{\prime}$ and $\mathcal{F}_{(A^{1},\ldots,A^{k})}^{%
\prime}$ (i.e., $*\equiv$\emph{Clifford product}) will be called
respectively the \emph{curl, scalar divergence, left contracted divergence }%
and\emph{\ gradient }of $\mathcal{F}_{(A^{1},\ldots,A^{k})}.$ Sometimes, $%
\mathcal{F}_{(A^{1},\ldots,A^{k})}^{\prime}$ will be called the \emph{%
standard derivative }of $\mathcal{F}_{(A^{1},\ldots,A^{k})}.$

We introduce now on the \emph{real vector space of differentially-induced }$%
r $\emph{-vector functionals of }$(p,q)$\emph{-extensor} the following \emph{%
four} derivative-like operators $\partial_{t}*$ as follows
\begin{equation}
\partial_{t}*\mathcal{F}_{(A^{1},\ldots,A^{k})}[t]=*\mathcal{F}%
_{(A^{1},\ldots,A^{k})}^{\prime}[t],   \label{7.4c}
\end{equation}
i.e., by eq.(\ref{7.4b})
\begin{equation}
\partial_{t}*\mathcal{F}_{(A^{1},\ldots,A^{k})}[t]=\overset{k}{\underset{i=1}%
{\sum}}A^{i}*\partial_{X^{i}}F[t(A^{1}),\ldots,t(A^{k})].   \label{7.4d}
\end{equation}

The special cases: $\partial_{t}\wedge,$ $\partial_{t}\cdot,$ $\partial
_{t}\lrcorner$ and $\partial_{t}$ (i.e., $*\equiv$\emph{Clifford product})
will be called respectively the (functional) \emph{curl, scalar divergence,
left contracted divergence} and \emph{gradient operator. }Sometimes, we will
say that $\partial_{t}$ is the \emph{standard derivative operator with
respect to }$t$\emph{.}

$\partial_{t}\wedge\mathcal{F}_{(A^{1},\ldots,A^{k})}[t],$ $\partial_{t}\cdot%
\mathcal{F}_{(A^{1},\ldots,A^{k})}[t],$ $\partial_{t}\lrcorner \mathcal{F}%
_{(A^{1},\ldots,A^{k})}[t]$ and $\partial_{t}\mathcal{F}_{(A^{1},%
\ldots,A^{k})}[t]$ (i.e., $*\equiv$\emph{Clifford product}) will be named
respectively as the \emph{curl, scalar divergence, left contracted
divergence }and \emph{gradient }of\emph{\ }$\mathcal{F}_{(A^{1},%
\ldots,A^{k})}.$ The gradient of $\mathcal{F}_{(A^{1},\ldots,A^{k})}$ will
be called the \emph{standard derivative }of $\mathcal{F}_{(A^{1},%
\ldots,A^{k})} $\emph{\ with respect to }$t$

It is still possible to define the noticeable derivative-like operator $%
A\cdot\partial_{t}$ as follows
\begin{align}
A\cdot\partial_{t}\mathcal{F}_{(A^{1},\ldots,A^{k})}[t] & =(A\cdot\frac {1}{%
p!}e^{j_{1}}\wedge\ldots e^{j_{p}})\mathcal{F}_{(A^{1},\ldots
,A^{k})e_{j_{1}}\wedge\ldots e_{j_{p}}}^{\prime}[t]  \label{7.4e} \\
& =(A\cdot\frac{1}{p!}e_{j_{1}}\wedge\ldots e_{j_{p}})\mathcal{F}%
_{(A^{1},\ldots,A^{k})e^{j_{1}}\wedge\ldots e^{j_{p}}}^{\prime}[t],
\label{7.4f}
\end{align}
i.e., by eq.(\ref{7.3c})
\begin{equation}
A\cdot\partial_{t}\mathcal{F}_{(A^{1},\ldots,A^{k})}[t]=\mathcal{F}%
_{(A^{1},\ldots,A^{k})A}^{\prime}[t].   \label{7.4g}
\end{equation}
Eq.(\ref{7.4g}) means that $A\cdot\partial_{t}$ is the $A$\emph{-directional
derivative operator} which maps $\mathcal{F}_{(A^{1},\ldots,A^{k})}$ $\mapsto
$ $\mathcal{F}_{(A^{1},\ldots,A^{k})A}^{\prime}.$

It is often convenient when doing calculations to employ some abuses of
notation for simplifying the handle of the fundamental formulas. Thus, eqs.(%
\ref{7.3a}) and (\ref{7.4d}) will be usually written
\begin{align}
A\cdot\partial_{t}F[t(A^{1}),\ldots,t(A^{k})] & =\overset{k}{\underset{i=1}{%
\sum}}A\cdot A^{i}\partial_{t(A^{i})}F[t(A^{1}),\ldots,t(A^{k})],
\label{7.4h} \\
\partial_{t}*F[t(A^{1}),\ldots,t(A^{k})] & =\overset{k}{\underset{i=1}{\sum }%
}A^{i}*\partial_{t(A^{i})}F[t(A^{1}),\ldots,t(A^{k})].   \label{7.4i}
\end{align}
No confusion arises since $A\cdot\partial_{t}$ and $\partial_{t}*$ denote
derivation of $r$-vector functional with respect to $(p,q)$-extensor $t,$
and $\partial_{t(A^{i})}$ holds for derivation of $r$-vector function with
respect to $q$-vector $t(A^{i})$.

It should be noted that by employing the abused notation we can re-write
eqs.(\ref{7.4aa}) and (\ref{7.4a}) as
\begin{align}
\partial_{t}*F[\ldots] & =\frac{1}{p!}(e^{j_{1}}\wedge\ldots
e^{j_{p}})*(e_{j_{1}}\wedge\ldots e_{j_{p}})\cdot\partial_{t}F[\ldots]
\label{7.4j} \\
& =\frac{1}{p!}(e_{j_{1}}\wedge\ldots e_{j_{p}})*(e^{j_{1}}\wedge\ldots
e^{j_{p}})\cdot\partial_{t}F[\ldots].   \label{7.4k}
\end{align}

\subsubsection{$A$-Directional Derivation Rules}

\begin{proposition}
Take a real $\lambda$ and a multivector $M.$ If $t\mapsto
F[t(A^{1}),\ldots,t(A^{k})]$ is any differentially-induced $r$-vector
functional of a $(p,q)$-extensor, then
\begin{align}
A\cdot\partial_{t}(\lambda F[\ldots]) & =\lambda A\cdot\partial_{t}F[\ldots],
\label{7.5a} \\
A\cdot\partial_{t}(F[\ldots]M) & =(A\cdot\partial_{t}F[\ldots])M.
\label{7.5b}
\end{align}
\end{proposition}

\begin{proof}
It follows directly from eq.(\ref{7.4h}) by using the derivation formulas:
$\partial_{X^{i}}(\lambda F(\ldots))=\lambda\partial_{X^{i}}F(\ldots)$ and
$\partial_{X^{i}}(F(\ldots)M)=(\partial_{X^{i}}F(\ldots))M.$
\end{proof}

\begin{theorem}
Let $t\mapsto F[t(A^{1}),\ldots,t(A^{k})]$ and $t\mapsto
G[t(A^{1}),\ldots,t(A^{k})]$ be any two differentially-induced $r$-vector
functionals of a $(p,q)$-extensor.

The addition $t\mapsto(F+G)[t(A^{1}),\ldots,t(A^{k})]$ is a
differentially-induced $r$-vector functional of a $(p,q)$-extensor and the
following rule holds
\begin{equation}
A\cdot\partial_{t}(F+G)[\ldots]=A\cdot\partial_{t}F[\ldots]+A\cdot\partial
_{t}G[\ldots].   \label{7.5c}
\end{equation}
\end{theorem}

\begin{proof}
As we can see, it is an immediate consequence of the derivation rule
$\partial_{X^{i}}(F+G)(\ldots)=\partial_{X^{i}}(F)(\ldots)+\partial_{X^{i}%
}G(\ldots).$
\end{proof}

\begin{theorem}
Let $t\mapsto\Phi[t(A^{1}),\ldots,t(A^{k})]$ and $t\mapsto
G[t(A^{1}),\ldots,t(A^{k})]$ be any differentially-induced scalar and $r$%
-vector functional of a $(p,q)$-extensor, respectively.

The scalar multiplication $t\mapsto(\Phi G)[t(A^{1}),\ldots,t(A^{k})]$ is
also a differentially-induced $r$-vector functional of a $(p,q)$-extensor
and we have
\begin{equation}
A\cdot\partial_{t}(\Phi G)[\ldots]=(A\cdot\partial_{t}\Phi[\ldots ]%
)G[\ldots]+\Phi[\ldots]A\cdot\partial_{t}G[\ldots].   \label{7.5d}
\end{equation}
\end{theorem}

It is rightly a Leibnitz-like rule.

\begin{proof}
As the reader can easily prove, eq.(\ref{7.5d}) is an immediate
consequence of
the derivation rule $\partial_{X^{i}}(\Phi G)(\ldots)=(\partial_{X^{i}}%
\Phi(\ldots))G(\ldots)+\Phi(\ldots)\partial_{X^{i}}G(\ldots).$
\end{proof}

\begin{theorem}
Let $t\mapsto\Psi[t(A^{1}),\ldots,t(A^{k})]$ and $\lambda\mapsto\phi(\lambda)
$ be any differentially-induced scalar functional and a derivable ordinary
real function, respectively. Then, $t\mapsto\phi(\Psi[t(A^{1}),%
\ldots,t(A^{k})])$ is a differentially-induced scalar functional and the
following rule holds
\begin{equation}
A\cdot\partial_{t}\phi(\Psi[\ldots])=\phi^{\prime}(\Psi[\ldots])A\cdot
\partial_{t}\Psi[\ldots].   \label{7.5e}
\end{equation}
\end{theorem}

It is an interesting and useful chain-like rule for $A$-directional
derivation of a special type of scalar functionals.

\begin{proof}
Eq.(\ref{7.5e}) follows easily from eq.(\ref{7.4h}) by taking into account the
derivation rule $\partial_{X^{i}}\phi\circ\Psi(\ldots)=\phi^{\prime}\circ
\Psi(\ldots)\partial_{X^{i}}\Psi(\ldots).$
\end{proof}

\section{Examples}

\begin{example}
{}Let $h\in ext_{1}^{1}(V)$ and take $a,b,c\in V$. Then,
\end{example}

\begin{align}
a\cdot\partial_{h}(h(b)\cdot h(c)) & =a\cdot b\partial_{h(b)}(h(b)\cdot
h(c))+a\cdot c\partial_{h(c)}(h(b)\cdot h(c))  \notag \\
& =a\cdot bh(c)+a\cdot ch(b),  \notag \\
a\cdot\partial_{h}(h(b)\cdot h(c)) & =h(a\cdot bc+a\cdot cb).   \label{7.6a}
\end{align}
Also,
\begin{align}
a\cdot\partial_{h}(h(b)\wedge h(c)) & =a\cdot b\partial_{h(b)}(h(b)\wedge
h(c))+a\cdot c\partial_{h(c)}(h(b)\wedge h(c))  \notag \\
& =a\cdot b(n-1)h(c)-a\cdot c(n-1)h(b)  \notag \\
& =(n-1)h(a\cdot bc-a\cdot cb),  \notag \\
a\cdot\partial_{h}(h(b)\wedge h(c)) & =(n-1)h(a\lrcorner(b\wedge c)).
\label{7.6b}
\end{align}
In eqs.(\ref{7.6a}) and (\ref{7.6b}) we have used the derivative formulas $%
\partial_{x}(x\cdot y)=y$ and $\partial_{x}(x\wedge y)=(n-1)y$, where $n$ is
the dimension of $V.$

The second formula developed in this example has an interesting and useful
generalization, which is:

The $a$-derivative of the $k$-vector functional $ext_{1}^{1}(V)\ni h\mapsto%
\underline{h}(a^{1}\wedge\ldots a^{k})\in\bigwedge^{k}V,$ with $%
a^{1},\ldots,a^{k}\in V,$ is given by
\begin{equation}
a\cdot\partial_{h}\underline{h}(a^{1}\wedge\ldots a^{k})=(n-k+1)\underline {h%
}(a\lrcorner(a^{1}\wedge\ldots a^{k})).   \label{7.6c}
\end{equation}

\begin{example}
Let $h\in ext_{1}^{1}(V)$ and take $b\in V.$

We shall calculate $a\cdot\partial_{h}h(b)$ and $a\cdot\partial_{h}h^{%
\dagger }(b).$ And, also $\partial_{h}*h(b)$ and $\partial_{h}*h^{%
\dagger}(b).$
\end{example}

First, we have
\begin{align}
a\cdot\partial_{h}h(b) & =a\cdot b\partial_{h(b)}h(b)=(a\cdot b)n,  \notag \\
a\cdot\partial_{h}h(b) & =n(a\cdot b),   \label{7.6d}
\end{align}
were we used the derivative formula $\partial_{x}x=n.$ Thus,
\begin{equation*}
\partial_{h}*h(b)=e^{j}*e_{j}\cdot\partial_{h}h(b)=e^{j}*n(e_{j}\cdot
b)=b*n,
\end{equation*}
i.e.,
\begin{align*}
\partial_{h}\wedge h(b) & =\partial_{h}h(b)=nb, \\
\partial_{h}\cdot h(b) & =\partial_{h}\lrcorner h(b)=0.
\end{align*}

Now, by employing a trick we have
\begin{equation*}
a\cdot\partial_{h}h^{\dagger}(b)=a\cdot\partial_{h}(h^{\dagger}(b)\cdot
e^{j}e_{j})=a\cdot\partial_{h}(b\cdot h(e^{j})e_{j}).
\end{equation*}

Thus, by using eq.(\ref{7.5b})
\begin{align}
a\cdot\partial_{h}h^{\dagger}(b) & =(\overset{n}{\underset{i=1}{\sum}}a\cdot
e^{i}\partial_{h(e^{i})}b\cdot h(e^{j}))e_{j}=\overset{n}{\underset{i=1}{%
\sum }}a\cdot e^{i}b\delta_{i}^{j}e_{j},  \notag \\
a\cdot\partial_{h}h^{\dagger}(b) & =ba,   \label{7.6e}
\end{align}
were we used the derivative formula $\partial_{x}(b\cdot x)=b.$ Thus,
\begin{equation*}
\partial_{h}*h^{\dagger}(b)=e^{j}*e_{j}\cdot\partial_{h}h^{%
\dagger}(b)=e^{j}*(be_{j}).
\end{equation*}
It follows that
\begin{align*}
\partial_{h}\wedge h^{\dagger}(b) & =e^{j}\wedge(b\cdot
e_{j})+e^{j}\wedge(b\wedge e_{j})=b. \\
\partial_{h}\cdot h^{\dagger}(b) & =e^{j}\cdot(b\cdot
e_{j})+e^{j}\cdot(b\wedge e_{j})=0. \\
\partial_{h}\lrcorner h^{\dagger}(b) & =e^{j}\lrcorner(b\cdot
e_{j})+e^{j}\lrcorner(b\wedge e_{j})=(e^{j}\cdot b)e_{j}-(e^{j}\cdot
e_{j})b=(1-n)b. \\
\partial_{h}h^{\dagger}(b) & =e^{j}(2e_{j}\cdot b-e_{j}b)=(2-n)b.
\end{align*}

\begin{example}
Let $t\in ext_{1}^{1}(V).$ The trace of $t,$ i.e., $t\mapsto
tr[t]=t(e^{j})\cdot e_{j},$ is a scalar functional and the bivector of $t,$
i.e., $t\mapsto biv[t]=t(e^{j})\wedge e_{j},$ is a bivector functional, both
of them associated to $t$. We shall calculate $a\cdot\partial_{t}tr[t]$ and $%
a\cdot\partial_{t}biv[t].$ And, also $\partial_{t}*tr[t]$ and $\partial
_{t}*biv[t].$
\end{example}

First, we have
\begin{align}
a\cdot\partial_{t}tr[t] & =\overset{n}{\underset{i=1}{\sum}}a\cdot
e^{i}\partial_{t(e^{i})}(t(e^{j})\cdot e_{j})=\overset{n}{\underset{i=1}{%
\sum }}a\cdot e^{i}\delta_{i}^{j}e_{j},  \notag \\
a\cdot\partial_{t}tr[t] & =a.   \label{7.6f}
\end{align}
We have used once again the derivative formula $\partial_{x}(x\cdot y)=y.$
Hence,
\begin{equation*}
\partial_{t}*tr[t]=e^{j}*e_{j}\cdot\partial_{t}tr[t]=e^{j}*e_{j},
\end{equation*}
i.e.,
\begin{align*}
\partial_{t}\wedge tr[t] & =0, \\
\partial_{t}\cdot tr[t] & =\partial_{t}\lrcorner tr[t]=\partial_{t}tr[t]=n.
\end{align*}

Now, we have also
\begin{align}
a\cdot\partial_{t}biv[t] & =\overset{n}{\underset{i=1}{\sum}}a\cdot
e^{i}\partial_{t(e^{i})}(t(e^{j})\wedge e_{j})=\overset{n}{\underset{i=1}{%
\sum}}a\cdot e^{i}(n-1)\delta_{i}^{j}e_{j},  \notag \\
a\cdot\partial_{t}biv[t] & =(n-1)a,   \label{7.6g}
\end{align}
were we have used once again the derivative formula $\partial_{x}(x\wedge
y)=(n-1)y.$ Hence,
\begin{equation*}
\partial_{t}*biv[t]=(n-1)e^{j}*e_{j},
\end{equation*}
i.e.,
\begin{align*}
\partial_{t}\wedge biv[t] & =0, \\
\partial_{t}\cdot biv[t] & =\partial_{t}\lrcorner biv[t]=\partial
_{t}biv[t]=(n-1)n.
\end{align*}

\begin{example}
Let $h\in ext_{1}^{1}(V)$ and take a non-zero $I\in\bigwedge^{n}V.$ We shall
calculate the $a$-directional derivative of the pseudoscalar functional $%
h\mapsto\underline{h}(I),$ i.e., $a\cdot\partial_{h}\underline{h}(I).$
\end{example}

By employing one of the expansion formulas for pseudoscalars (see \cite{4}),
eq.(\ref{7.5a}) and eq.(\ref{7.6c}) we have
\begin{align}
a\cdot\partial_{h}\underline{h}(I) &
=a\cdot\partial_{h}I\cdot(e_{1}\wedge\ldots e_{n})\underline{h}%
(e^{1}\wedge\ldots e^{n})  \notag \\
& =I\cdot(e_{1}\wedge\ldots e_{n})a\cdot\partial_{h}\underline{h}%
(e^{1}\wedge\ldots e^{n})  \notag \\
& =I\cdot(e_{1}\wedge\ldots e_{n})\underline{h}(a\lrcorner(e^{1}\wedge\ldots
e^{n})),  \notag \\
a\cdot\partial_{h}\underline{h}(I) & =\underline{h}(a\lrcorner I)=h(aI).
\label{7.6h}
\end{align}

\begin{example}
Let $h\in ext_{1}^{1}(V)$ and take a non-zero $I\in\Lambda^{n}V.$ The
determinant of $h,$ i.e., $h\mapsto\det[h]$ such that $\underline{h}(I)=\det[%
h]I,$ is a characteristic scalar functional of $h.$ We shall calculate $%
a\cdot\partial_{h}\det[h]$ and $\partial_{h}*\det[h].$
\end{example}

By employing eq.(\ref{7.5b}) and eq.(\ref{7.6h}) we have
\begin{equation*}
a\cdot\partial_{h}\det[h]=(a\cdot\partial_{h}\underline{h}(I))I^{-1}=%
\underline{h}(aI)I^{-1}.
\end{equation*}
But, by taking into account the extensor formula $h^{-1}(a)=\left.
\det\right. ^{-1}[h]\underline{h}^{\dagger}(aI)I^{-1}$ (see\cite{5}) and
recalling that $\det[h^{\dagger}]=\det[h]$ and $h^{*}=(h^{%
\dagger})^{-1}=(h^{-1})^{\dagger}$ we get
\begin{equation}
a\cdot\partial_{h}\det[h]=\det[h]h^{*}(a).   \label{7.6i}
\end{equation}

Hence, it follows that
\begin{equation*}
\partial_{h}*\det[h]=e^{j}*e_{j}\cdot\partial_{h}\det[h]=\det[h]%
e^{j}*h^{*}(e_{j}),
\end{equation*}
i.e.,
\begin{align*}
\partial_{h}\wedge\det[h] & =-\det[h]h^{*}(e_{j})\wedge e^{j}=\det
[h]biv[h^{-1}]. \\
\partial_{h}\cdot\det[h] & =\partial_{h}\lrcorner\det[h]=\det[h]%
h^{-1}(e^{j})\cdot e_{j}=\det[h]tr[h^{-1}]. \\
\partial_{h}\det[h] & =\det[h]e^{j}h^{*}(e_{j})=\det[h]%
(tr[h^{-1}]+biv[h^{-1}]).
\end{align*}

\subsection{An Enlightening Discussion}

Let us consider for example a differentially-induced scalar functional of $%
(1,1)$-extensor $t\mapsto\Phi[t(a^{1})].$ We have the possibility for
constructing a differentiable scalar function of $n\times n$ real variables $%
(t_{11},\ldots,t_{1n},\ldots,t_{n1},\ldots,t_{nn})\mapsto\widehat{\Phi }%
(t_{11},\ldots,t_{1n},\ldots,t_{n1},\ldots,t_{nn}),$ defined by
\begin{equation}
\widehat{\Phi}(t_{11},\ldots,t_{1n},\ldots,t_{n1},\ldots,t_{nn})=\Phi
[t_{ij}(a^{1}\cdot e^{i})e^{j}]   \label{7.7a}
\end{equation}
where $t_{ij}=t(e_{i})\cdot e_{j}$ are the $n\times n$ matrix elements of $t$
with respect to $\{e_{k}\}.$

Eq.(\ref{7.7a}) shows that all information just contained into the \emph{%
classical real function} $(t_{11},\ldots,t_{nn})\mapsto\widehat{\Phi }%
(t_{11},\ldots,t_{nn})$ whose real variables are $t_{pq},$ is also codified
into the scalar functional $t\mapsto\Phi[t(a^{1})].$

We shall search for the relationship which exists between the \emph{ordinary
partial derivatives} of $\widehat{\Phi}(\ldots)$ with respect to each \emph{%
tensor covariant component}\footnote{%
They are the $n\times n$ covariant components of a $2$-tensor $T$ in
biunivocal correspondence with the $(1,1)$-extensor $t,$ see \cite{5}, i.e.,
$T_{pq}\equiv T(e_{p},e_{q})=t(e_{p})\cdot e_{q}\equiv t_{pq}.$} $t_{pq}$
and the $a$-directional derivative of $\Phi[\ldots].$

By using $\partial_{\lambda_{i}}\Phi(x(\lambda_{1},\ldots,\lambda
_{k}))=\partial_{\lambda_{i}}x(\lambda_{1},\ldots,\lambda_{k})\cdot
\partial_{x}\Phi(x(\lambda_{1},\ldots,\lambda_{k})),$ a chain-like
derivation rule, we may write
\begin{align}
\dfrac{\partial\widehat{\Phi}}{\partial t_{pq}}(t_{11},\ldots,t_{nn}) &
=\partial_{t_{pq}}(t_{ij}(a^{1}\cdot e^{i})e^{j})\cdot\partial_{x}\Phi
[t_{ij}(a^{1}\cdot e^{i})e^{j}]  \notag \\
& =\delta_{ij}^{pq}(a^{1}\cdot e^{i})e^{j}\cdot\partial_{t(a^{1})}\Phi[%
t(a^{1})],  \notag \\
\dfrac{\partial\widehat{\Phi}}{\partial t_{pq}}(t_{11},\ldots,t_{nn}) &
=(a^{1}\cdot e^{p})e^{q}\cdot\partial_{t(a^{1})}\Phi[t(a^{1})].
\label{7.7b}
\end{align}
Now, Clifford multiplication by $(a\cdot e_{p})e_{q}$ (and summing over $p,q$%
) on both sides of eq.(\ref{7.7b}) yields
\begin{align}
(a\cdot e_{p})e_{q}\dfrac{\partial\widehat{\Phi}}{\partial t_{pq}}%
(t_{11},\ldots,t_{nn}) & =(a\cdot a^{1})e_{q}e^{q}\cdot\partial_{t(a^{1})}%
\Phi[t(a^{1})]  \notag \\
& =a\cdot a^{1}\partial_{t(a^{1})}\Phi[t(a^{1})],  \notag \\
(a\cdot e_{p})e_{q}\dfrac{\partial\widehat{\Phi}}{\partial t_{pq}}%
(t_{11},\ldots,t_{nn}) & =a\cdot\partial_{t}\Phi[t(a^{1})].   \label{7.7c}
\end{align}
That is the required result relating both $\dfrac{\partial\widehat{\Phi}}{%
\partial t_{pq}}(t_{11},\ldots,t_{nn})$ and $a\cdot\partial_{t}\Phi
[t(a^{1})].$

It is still possible to find a relationship between $\dfrac{\partial
\widehat{\Phi}}{\partial t_{pq}}(t_{11},\ldots,t_{nn})$ and the $*$%
-derivatives of $\Phi[t(a^{1})].$ From eq.(\ref{7.7b}) we have
\begin{align}
e_{p}*(e_{q}\dfrac{\partial\widehat{\Phi}}{\partial t_{pq}}%
(t_{11},\ldots,t_{nn})) & =a^{1}*(e_{q}e^{q}\cdot\partial_{t(a^{1})}\Phi
[t(a^{1})])  \notag \\
& =a^{1}*\partial_{t(a^{1})}\Phi[t(a^{1})],  \notag \\
e_{p}*(e_{q}\dfrac{\partial\widehat{\Phi}}{\partial t_{pq}}%
(t_{11},\ldots,t_{nn})) & =\partial_{t}*\Phi[t(a^{1})].   \label{7.7d}
\end{align}

That is the expected identity which relates both $\dfrac{\partial\widehat {%
\Phi}}{\partial t_{pq}}(t_{11},\ldots,t_{nn})$ and $\partial_{t}*\Phi
[t(a^{1})].$

\section{Conclusions}

In this paper we introduced the key concepts of a theory of multivector
functionals. We studied several aspects of the notion of derivative that can
be applied to these objects, as e.g., the $A$-directional derivatives and
the generalized concepts of curl, divergence and gradient. We worked in
details several examples where we calculate different types of derivatives
for multivector functionals. It is worth to said once again that these
objects play a decisive role in the development of a Lagrangian formalism
for extensor fields as it will be seen in two future series of papers: \emph{%
geometric theories of gravitation and Lagrangian formulation of the
multivector and extensor fields theory}.\bigskip

\textbf{Acknowledgement}: V. V. Fern\'{a}ndez is grateful to FAPESP for a
posdoctoral fellowship. W. A. Rodrigues Jr. is grateful to CNPq for a senior
research fellowship (contract 201560/82-8) and to the Department of
Mathematics of the University of Liverpool for the hospitality. Authors are
also grateful to Drs. P. Lounesto, I. Porteous, and J. Vaz, Jr. for their
interest on our research and useful discussions.

\end{document}